\documentclass[11pt,reqno]{amsart}

\usepackage{amsmath, amsfonts, amsthm, amssymb}

\newtheorem{Theorem}{Theorem}[section]
\newtheorem{Lemma}{Lemma}[section]

\theoremstyle{definition}

\theoremstyle{remark}
\newtheorem{Remark}{Remark}[section]

\numberwithin{equation}{section}

\renewcommand{\r}{\rho}

\renewcommand{\u}{{\bf u}}

\newcommand{\R}{{\mathbb R}}
\newcommand{\Dv}{{\rm div}}
\newcommand{\tr}{{\rm tr}}

\newcommand{\dl}{\delta}

\def\f{\frac}

\def\ov{\overline}

\def\D{\Delta }

\def\hf1{^\f{1}{1-\xi^2}}

\def\be{\begin{equation}}
\def\en{\end{equation}}
\def\bs{\begin{split}}
\def\es{\end{split}}

\newcommand{\F}{{\mathtt F}}

\author{Xianpeng HU and Dehua Wang}
\address{Department of Mathematics, University of Pittsburgh,
                           Pittsburgh, PA 15260, USA.}
\email{xih15@pitt.edu}
\address{Department of Mathematics, University of Pittsburgh,
                           Pittsburgh, PA 15260, USA.}
\email{dwang@math.pitt.edu}

\title[Local Strong Solution to Compressible Viscoelastic Fluid]
{Local Strong Solution to the Compressible Viscoelastic Fluid
with Large Data}

\date{\today}

\begin{document}

\begin{abstract}
The existence and uniqueness of local in time strong solution with
large initial data for the three-dimensional compressible
viscoelastic fluid is established. The strong solution has weaker
regularity than the classical solution. The Lax-Milgram theorem
and the Schauder-Tychonoff fixed-point argument are applied.
\end{abstract}
\maketitle

\section{Introduction}
Elastic solids and viscous fluids are two extremes of material behavior.
Viscoelastic fluids  show intermediate behavior with  some remarkable phenomena  due to their ``elastic" nature.
These fluids exhibit a combination of both fluid and solid characteristics, keep memory of their past deformations,
and their behaviour is a function of these old deformations.
Viscoelastic fluids have a wide range of applications  and hence have received a great deal of interest.
Examples and applications of viscoelastic fluids include from oil, liquid polymers,
mucus, liquid soap,  toothpaste,  clay,   ceramics,  gels,  some types of suspensions,
to  bioactive fluids,  coatings and drug delivery systems for controlled drug release, scaffolds for tissue engineering,
and viscoelastic blood fluid flow past valves; see \cite{GBO, KM, YSS} for more applications.
For the viscoelastic materials, the competition between the kinetic energy and the internal elastic
energy through the special transport properties of their respective internal elastic variables makes the materials more
untractable in  understanding their behavior, since any distortion of
microstructures, patterns or configurations in the dynamical flow
will involve the deformation tensor. For classical simple fluids,  the internal energy can be
determined solely by the determinant of the deformation tensor; however, the internal energy of complex fluids carries all the
information of the deformation tensor. The interaction between the microscopic elastic properties and the macroscopic fluid
motions  leads to the rich and complicated rheological phenomena in viscoelastic fluids, and also causes
formidable analytic and numerical challenges in mathematical analysis.
The equations of the compressible viscoelastic fluids of Oldroyd type (\cite{Oldroyd1, Oldroyd2})
in three spatial dimensions take the following form \cite{Joseph, LW, RHN}:
\begin{subequations} \label{e1}
\begin{align}
&\r_t +\Dv(\r\u)=0,\label{e11}\\
&(\r\u)_t+\Dv\left(\r\u\otimes\u\right)-\mu\D \u-(\lambda+\mu)\nabla\Dv\u+\nabla P(\r)=\Dv(\r \F\F^\top),\label{e12}\\
&\F_t+\u\cdot\nabla\F=\nabla\u\, \F,\label{e13}
\end{align}
\end{subequations}
where $\r$ stands for the density, $\u\in \R^3$ the velocity,
$\F\in M^{3\times 3}$ the deformation gradient,  and $P(\r)$ the
pressure which is a strictly increasing convex function of the
density. The notation $M^{3\times 3}$ means the set of all
$3\times 3$ matrixes. The viscosity coefficients $\mu$ and
$\lambda$ satisfy the conditions that $3\mu+2\lambda>0$ and
$\mu>0$, which ensure that the operator $-\mu\D
\u-(\lambda+\mu)\nabla\Dv\u$ is a strongly elliptic operator. The
symbol $\otimes$ denotes the Kronecker tensor product,  $\F^\top$
means the transpose matrix of $\F$, and the notation
$\u\cdot\nabla\F$ is understood to be $(\u\cdot\nabla)\F$.
Usually, we refer the equation \eqref{e11} as the continuity
equation.

We are interested in the Cauchy problem of \eqref{e1} with the initial condition:
\begin{equation}\label{IC}
(\r, \u, \F)(0,x)=(\r_0, \u_0, \F_0)(x), \quad x\in\R^3.
\end{equation}
The aim of  this paper is to establish the local existence and
uniqueness of strong solution to   system \eqref{e1} with large
initial data in the three dimensional space $\R^3$. By a
$\textit{strong solution}$, we mean a triplet $(\r, \u, \F)$ with
$\u(t,\cdot)\in W^{2,q}$ and $(\r(t,\cdot), \F(t,\cdot))\in
W^{1,q}, \, q>3$ satisfying \eqref{e1} almost everywhere
with the initial condition \eqref{IC}. There have been many
studies and rich results in the literature for the global
existence of classical solutions (namely in $H^3$ or other
functional spaces with much higher regularity) for the
corresponding {\em incompressible} viscoelastic fluids, see
 \cite{CM, CZ, KP, LLZ, LLZH2, LLZH, LZP, LW, LZ} and the references therein.
 For the compressible viscoelastic fluids \eqref{e1},  the global existence of classical solutions in $H^3$
with small perturbation near its equilibrium for \eqref{e1}
without the pressure term was studied in \cite{LLZH31}, and a
local existence of strong solution near the equilibrium and a
series of uniform estimates were obtained in \cite{HW}. One of the
main difficulties in proving the global existence is the lacking
of the dissipative estimates for the deformation gradient and the
gradient of the density. To overcome this difficulty, for
incompressible cases, authors in \cite{LLZH} introduced an
auxiliary function to obtain the dissipative estimate for the
classical solutions, while authors in \cite{LZ} directly deal with
the quantities such as $\Delta \u+\Dv \F$.
Since we are concerned with the {\em strong solutions in
$W^{2,q}$} which have weaker regularity than the {\em classical
solution in $H^3$}, we first linearize \eqref{e1} and use the
Lax-Milgram theorem to obtain the solution to the linearized system, then
we apply the Schauder-Tychonoff fixed point theorem to obtain the
strong solution of \eqref{e1}.
%
Using the standard notations  $W^{s, q}(\R^3)$ ($H^s(\R^3)$ if $q=2$) for the Sobolev
spaces and setting
$Q_T:=[0,T]\times\R^3$ for any $T>0$, we can state our result on the existence and uniqueness as follows.

\begin{Theorem}\label{T20}
Assume that
$$\r_0\in W^{1,q}(\R^3)\cap H^1(\R^3),\quad \u_0\in H^2(\R^3),\quad  \F_0\in W^{1,q}(\R^3)\cap H^1(\R^3),$$
for some  $q\in(3, 6]$, and for some positive constants $\alpha, \beta$, and $r_0$,
$$\alpha\le \r_0\le\beta,\quad
\|\u_0\|_{H^2}+\|\r_0\|_{W^{1,q}\cap H^1}+\|\F_0\|_{W^{1,q}\cap H^1} \le r_0.$$
 Then there are positive constants $\overline{T}=\ov{T}(r_0)$, $\alpha (\ov{T}, r_0,
\alpha)$, and $\beta(\ov{T}, r_0, \beta)$, such that,  the Cauchy problem \eqref{e1}-\eqref{IC} has a unique strong solution
$(\r,\u, \F)$ defined for $(t,x)\in(0,\ov{T})\times\R^3$,  satisfying
$$\alpha(\ov{T},r_0, \alpha)\le\r\le\beta(\ov{T},r_0,\beta);$$
$$\r\in L^\infty(0,\ov{T}; W^{1,q}(\R^3)\cap H^1(\R^3))\cap
L^\infty(Q_{\ov{T}});\quad \partial_t\r\in L^\infty(0,\ov{T};
L^q(\R^3));$$
$$\u\in L^2(0,\ov{T}; W^{2,q}(\R^3)\cap H^2(\R^3)); \quad \partial_t\u\in
L^2(0,\ov{T}; H^1(\R^3));$$
$$\F\in L^\infty(0,\ov{T}; W^{1,q}(\R^3)\cap H^1(\R^3))\cap
L^\infty(Q_{\ov{T}});\quad \partial_t\F\in L^\infty(0,\ov{T};
L^q(\R^3)).$$
\end{Theorem}

 The viscoelasticity system \eqref{e1}
can be regarded as a combination between the compressible
Navier-Stokes equation with the source term $\Dv(\r\F\F^\top)$ and
the equation \eqref{e13}. As for the global existence of classical
solutions of the small perturbation near an equilibrium for
compressible Navier-Stokes equations, we refer the interested
reader to \cite{MT1, MT} and the references cited therein. The
global existence of strong solutions with small perturbations near
an equilibrium for compressible Navier-Stokes equations was also
discussed in \cite{AI, SS}. Unlike the  Navier-Stokes equations (\cite{f2,p1,p2}),
the global existence of weak solutions to \eqref{e1} with large
initial data is still an outstanding open question. In this
direction for the incompressible viscoelasticity, when the
contribution of the strain rate (symmetric part of $\nabla\u$) in
the constitutive equation is neglected, Lions and Masmoudi in
\cite{LM} proved the global existence of weak solutions with large
initial data for the Oldroyd model. Also Lin, Liu, and Zhang
showed in \cite{LLZ} the existence of global weak solutions with
large initial data for the incompressible viscoelasticity if the
velocity satisfies the Lipschitz condition. When dealing with the
global existence of weak solutions with large data in the
compressible case, among all of difficulties, the rapid
oscillation of the density and the non-compatibility between the
quadratic form and the weak convergence are of the main issues.
For the inviscid elastodynamics, see \cite{FJ,ST,ST2} and their references on the finite-time blow up and global existence of classical solutions.

This paper will be organized as follows. In Section 2,  we will
recall briefly the compressible viscoelastic system from some
fundamental mechanical theory. In Section 3 and Section 4, we will
give the proof of the main theorem (Theorem \ref{T20}). More precisely, Section 3 is
devoted to the local existence of the system \eqref{e1} by the
Lax-Milgram theorem and the Schauder-Tychonoff fixed-point
argument, while the main goal of Section 4 is to prove the
uniqueness of the solution obtained in Section 3.

\bigskip\bigskip

\section{Mechanical Background of Viscoelasticity}

In this section, we discuss some background of viscoelasticity and some basic properties of system \eqref{e1}.
First, we recall the definition of the deformation
gradient $\F$. The dynamics of any mechanical problem with a
velocity field $\u(x,t)$ can be described by the flow map, a time
dependent family of orientation preserving diffeomorphisms
$x(t,X)$, $0\le t\le T$:
\begin{equation}\label{x2}
\begin{cases}
&\f{d}{dt}x(t,X)=\u(t, x(t,X)),\\
&x(0,X)=X.
\end{cases}
\end{equation}
The material point (labelling) $X$ in the reference configuration
is deformed to the spatial position $x(t,X)$ at time $t$, which is
the observer's coordinate.

The deformation gradient $\widetilde{\F}$ is used to describe the
changing of any configuration, amplification or pattern during the
dynamical process, which is defined as
$$\widetilde{\F}(t,X)=\f{\partial x}{\partial X}(t,X).$$
Notice that this quantity is defined in the Lagrangian material
coordinate. Obviously it satisfies the following rule, by changing
the order of the differentiation:
\begin{equation}\label{x1}
\f{\partial\widetilde{\F}(t,X)}{\partial
t}=\f{\partial\u(t,x(t,X))}{\partial X}.
\end{equation}
In the Eulerian coordinate, the corresponding deformation gradient
$\F(t,x)$ will be defined as $F(t, x(t,X))=\widetilde{\F}(t,X)$.
The equation \eqref{x1}, together with the chain rule and
\eqref{x2}, yield the following equation:
\begin{equation*}
\begin{split}
\partial_t\F(t, x(t,X))+\u\cdot\nabla\F(t,x(t,X))&=\partial_t\F(t, x(t,X))+\f{\partial\F(t,x(t,X))}{\partial x}
\cdot\f{\partial x(t,X)}{\partial t}\\
&=\f{\partial\widetilde{\F}(t,X)}{\partial t}=\f{\partial\u(t,x(t,X))}{\partial X}\\
&=\f{\partial\u(t, x(t,X))}{\partial x}\f{\partial x}{\partial X}\\
&=\f{\partial\u(t, x(t,X))}{\partial x}\widetilde{\F}(t,
X)=\nabla\u\cdot\F,
\end{split}
\end{equation*}
which is exactly the equation \eqref{e13}. Here, and in what follows,
we use the conventional notations:
$$(\nabla\u)_{ij}=\f{\partial u_i}{\partial x_j},\quad(\nabla\u\;\F)_{i,j}=(\nabla\u)_{ik}F_{kj},\quad (\u\cdot\nabla\F)_{ij}
=u_k\f{\partial \F_{ij}}{\partial x_{k}},$$ and summation over
repeated indices will always be well understood. In
viscoelasticity, \eqref{e13} can also be interpreted as the
consistency of the flow maps generated by the velocity field $\u$
and the deformation gradient $\F$.

The difference between fluids and solids lies in the fact that in
fluids, such as Navier-Stokes equations \cite{AI}, the internal
energy can be determined solely by the determinant part of $\F$
(equivalently the density $\r$, and hence, \eqref{e13} can be
disregarded) and in elasticity, the energy depends on all
information of $\F$.

In the continuum physics, if we assume that the material is
homogeneous, the conservation laws of mass and of momentum become
\cite{CD, LLZH, ST2}
\begin{equation}\label{x3}
\partial_t\r+\Dv(\r\u)=0,
\end{equation}
and
\begin{equation}\label{x4}
\partial_t(\r\u)+\Dv(\r\u\otimes\u)-\mu\Delta\u-(\lambda+\mu)\nabla\Dv\u+\nabla P(\r)=\Dv((\det\F)^{-1}S\F^\top),
\end{equation}
where
\begin{equation}\label{x5}
\r\det\F=1,
\end{equation}
and
\begin{equation}\label{y1}
S_{ij}(\F)=\f{\partial W}{\partial F_{ij}}.
\end{equation}
Here $S$, $\r S\F^\top$, $W(\F)$ denote $\textit{Piola-Kirchhoff
stress}$, $\textit{Cauchy stress}$ and the elastic energy of the
material respectively. Recall that the condition \eqref{y1}
implies that the material is called hyperelastic \cite{LW}. In the
case of Hookean (linear) elasticity \cite{LLZH2, LLZH, LZP},
\begin{equation}\label{x6}
W(\F)=\f{1}{2}|\F|^2=\f{1}{2}tr(\F\F^\top),
\end{equation}
where the notation $\tr$ stands for the trace operator of a
matrix, and hence,
\begin{equation}\label{x7}
S(\F)=\F.
\end{equation}
Combining  equations \eqref{x2}-\eqref{x7} together, we
obtain system \eqref{e1}.

If the viscoelastic system \eqref{e1} satisfies
$\Dv(\r_0\F_0^\top)=0$, it is verified in \cite{LZ} (see
Proposition 3.1) that this condition will insist in time, that is,
\begin{equation}\label{cc1}
\Dv(\r(t)\F(t)^\top)=0,\quad\textrm{for}\quad t\ge 0.
\end{equation}

Another hidden, but important, property of the viscoelastic fluids
\eqref{e1} is concerned with the curl of the deformation gradient
(for the incompressible case, see \cite{LLZH2, LLZH}). Actually,
the following lemma says that the curl of the deformation gradient
is of higher order.

\begin{Lemma}\label{curl}
Assume that \eqref{e13} is satisfied and $(\u, \F)$ is the
solution of the system \eqref{e1}. Then the following identity
\begin{equation}\label{curl1}
\F_{lk}\nabla_l \F_{ij}=\F_{lj}\nabla_l \F_{ik}
\end{equation}
holds for all time $t> 0$ if it initially satisfies \eqref{curl1}.
\end{Lemma}

Again, here, the standard summation notation over the repeated
index is adopted.

\begin{proof}
First, we establish the evolution equation for the equality
$$\F_{lk}\nabla_l \F_{ij}-\F_{lj}\nabla_l \F_{ik}.$$
 Indeed, by the equation \eqref{e13}, we have
\begin{equation*}
\partial_t\nabla_l\F_{ij}+\u\cdot\nabla\nabla_l
\F_{ij}+\nabla_l\u\cdot\nabla\F_{ij}=\nabla_m\u_i\nabla_l\F_{mj}+\nabla_l\nabla_m\u_i
\F_{mj}.
\end{equation*}
Thus,
\begin{equation}\label{curl2}
\F_{lk}(\partial_t\nabla_l\F_{ij}+\u\cdot\nabla\nabla_l
\F_{ij})+\F_{lk}\nabla_l\u\cdot\nabla\F_{ij}=\F_{lk}\nabla_m\u_i\nabla_l\F_{mj}+\F_{lk}\nabla_l\nabla_m\u_i
\F_{mj}.
\end{equation}
Also, from \eqref{e13}, we obtain
\begin{equation}\label{curl3}
\nabla_l\F_{ij}(\partial_t\F_{lk}+\u\cdot\nabla\F_{lk})=\nabla_l\F_{ij}\nabla_m\u_l\F_{mk}.
\end{equation}

Now, adding \eqref{curl2} and \eqref{curl3}, we deduce that
\begin{equation}\label{curl4}
\begin{split}
\partial_t(\F_{lk}\nabla_l \F_{ij})+\u\cdot\nabla(\F_{lk}\nabla_l
\F_{ij})&=-\F_{lk}\nabla_l\u\cdot\nabla\F_{ij}+\F_{lk}\nabla_m\u_i\nabla_l\F_{mj}\\&\quad+\F_{lk}\nabla_l\nabla_m\u_i
\F_{mj}+\nabla_l\F_{ij}\nabla_m\u_l\F_{mk}\\
&=\F_{lk}\nabla_m\u_i\nabla_l\F_{mj}+\F_{lk}\nabla_l\nabla_m\u_i
\F_{mj}.
\end{split}
\end{equation}
Here, we used the identity which is derived by interchanging the
roles of indices $l$ and $m$:
$$\F_{lk}\nabla_l\u\cdot\nabla\F_{ij}=\F_{lk}\nabla_l\u_m\nabla_m\F_{ij}=\nabla_l\F_{ij}\nabla_m\u_l\F_{mk}.$$
Similarly, one has
\begin{equation}\label{curl5}
\begin{split}
\partial_t(\F_{lj}\nabla_l \F_{ik})+\u\cdot\nabla(\F_{lj}\nabla_l
\F_{ik})=\F_{lj}\nabla_m\u_i\nabla_l\F_{mk}+\F_{lj}\nabla_l\nabla_m\u_i
\F_{mk}.
\end{split}
\end{equation}
Subtracting \eqref{curl5} from \eqref{curl4} yields
\begin{equation}\label{curl6}
\begin{split}
&\partial_t(\F_{lk}\nabla_l \F_{ij}-\F_{lj}\nabla_l
\F_{ik})+\u\cdot\nabla(\F_{lk}\nabla_l \F_{ij}-\F_{lj}\nabla_l
\F_{ik})\\&\quad=\nabla_m\u_i(\F_{lk}\nabla_l\F_{mj}-\F_{lj}\nabla_l\F_{mk})+\nabla_l\nabla_m\u_i
(\F_{mj}\F_{lk}-\F_{mk}\F_{lj}).
\end{split}
\end{equation}
Due to the fact
$$\nabla_l\nabla_m\u_i=\nabla_m\nabla_l\u_i$$ in the sense of distributions, we
have, again by interchanging the roles of indices $l$ and $m$,
\begin{equation*}
\begin{split}
\nabla_l\nabla_m\u_i
(\F_{mj}\F_{lk}-\F_{mk}\F_{lj})&=\nabla_l\nabla_m\u_i
\F_{mj}\F_{lk}-\nabla_l\nabla_m\u_i
\F_{mk}\F_{lj}\\
&=\nabla_l\nabla_m\u_i \F_{mj}\F_{lk}-\nabla_m\nabla_l\u_i
\F_{lk}\F_{mj}\\
&=(\nabla_l\nabla_m\u_i-\nabla_m\nabla_l\u_i) \F_{lk}\F_{mj}=0.
\end{split}
\end{equation*}
From this identity, equation \eqref{curl6} can be simplified as
\begin{equation}\label{curl7}
\begin{split}
&\partial_t(\F_{lk}\nabla_l \F_{ij}-\F_{lj}\nabla_l
\F_{ik})+\u\cdot\nabla(\F_{lk}\nabla_l \F_{ij}-\F_{lj}\nabla_l
\F_{ik})\\&\quad=\nabla_m\u_i(\F_{lk}\nabla_l\F_{mj}-\F_{lj}\nabla_l\F_{mk}).
\end{split}
\end{equation}
Multiplying \eqref{curl7} by $\F_{lk}\nabla_l
\F_{ij}-\F_{lj}\nabla_l \F_{ik}$, we get
\begin{equation}\label{curl8}
\begin{split}
&\partial_t|\F_{lk}\nabla_l \F_{ij}-\F_{lj}\nabla_l
\F_{ik}|^2+\u\cdot\nabla|\F_{lk}\nabla_l \F_{ij}-\F_{lj}\nabla_l
\F_{ik}|^2\\&\quad=2(\F_{lk}\nabla_l \F_{ij}-\F_{lj}\nabla_l
\F_{ik})\nabla_m\u_i(\F_{lk}\nabla_l\F_{mj}-\F_{lj}\nabla_l\F_{mk})\\
&\quad\le 2\|\nabla\u\|_{L^\infty(\R^3)}\mathcal{M}^2,
\end{split}
\end{equation}
where $\mathcal{M}$ is defined as
$$\mathcal{M}=\max_{i,j,k}\{|\F_{lk}\nabla_l
\F_{ij}-\F_{lj}\nabla_l \F_{ik}|^2\}.$$ Hence, \eqref{curl8}
implies
\begin{equation}\label{curl9}
\partial_t\mathcal{M}+\u\cdot\nabla\mathcal{M}\le
2\|\nabla\u\|_{L^\infty(\R^3)}\mathcal{M}.
\end{equation}

On the other hand, the characteristics of $\partial_t
f+\u\cdot\nabla f=0$ is given by
$$\f{d}{ds}X(s)=\u(s,X(s)),\quad X(t)=x.$$
Hence, \eqref{curl8} can be rewritten as
\begin{equation}\label{curl10}
\f{\partial U}{\partial t}\le B(t,y)U,\quad
U(0,y)=\mathcal{M}_0(y),
\end{equation}
where
$$U(t,y)=\mathcal{M}(t, X(t,x)),\quad
B(t,y)=2\|\nabla\u\|_{L^\infty(\R^3)}(t, X(t,y)).$$ The
differential inequality \eqref{curl10} implies that
$$U(t,y)\le U(0)\exp\left(\int_0^t B(s,y)ds\right).$$ Hence,
$$\mathcal{M}(t,x)\le \mathcal{M}(0)\exp\left(\int_0^t
2\|\nabla\u\|_{L^\infty(\R^3)}(s)ds\right).$$ Hence, if
$\mathcal{M}(0)=0$, then $\mathcal{M}(t)=0$ for all $t>0$, and
 the proof of the lemma is complete.
\end{proof}

\begin{Remark}
Lemma \ref{curl} can be interpreted from the physical viewpoint as
follows: formally, the fact that the Lagrangian derivatives
commute and the definition of the deformation gradient imply
$$\partial_{X_k}\widetilde{\F}_{ij}=\f{\partial^2x_i}{\partial X_k\partial X_j}=\f{\partial^2x_i}{\partial X_j\partial X_k}
=\partial_{X_j}\widetilde{\F}_{ik},$$ which is equivalent to, in
the Eulerian coordinates,
$$\widetilde{\F}_{lk}\nabla_{l}\F_{ij}(t,x(t,X))=\widetilde{\F}_{lj}\nabla_{l}\F_{ik}(t, x(t,X)), $$
that is,
$$\F_{lk}\nabla_{l}\F_{ij}(t,x)=\F_{lj}\nabla_{l}\F_{ik}(t,x).$$
\end{Remark}

Finally, if the density $\r$ is a constant, the equations of the
incompressible viscoelastic fluids have the following form (see
\cite{CZ, LLZ, LLZH2, LLZH, LZP, LZ} and references therein):
\begin{equation}\label{ICV}
\begin{cases}
\Dv \u=0,\\
\partial_t\u+\u\cdot\nabla\u-\mu\Delta\u+\nabla P=\Dv(\F\F^\top),\\
\partial_t\F+\u\cdot\nabla\F=\nabla\u\,\F.
\end{cases}
\end{equation}

\bigskip\bigskip

\section{Local Existence}
In this section, we will prove the existence part in Theorem
\ref{T20}. The proof will proceed  through four steps by combining
the Lax-Milgram theorem and a fixed point argument. To this end,
we consider first the linearized problem.

Set
\begin{equation*}
\begin{split}
\mathcal{W}'=\big\{ \psi\in (L^2(0,T; H^2(\R^3)))^3: \;
\partial_t\psi\in L^2(Q_T)\big\}
\end{split}
\end{equation*}
with the natural norm $\|\psi\|_{\mathcal{W}'}$,  and for $q\in
(3,6]$ we define
\begin{equation*}
\begin{split}
\mathcal{W}=\mathcal{W}'&\cap \left(L^2(0,T; W^{2,q}(\R^3))\cap L^\infty(0,T; H^2(\R^3))\right)^3\\
&\cap\{\psi\in (L^2(0,T; H^2(\R^3)))^3: \; \psi_t\in L^\infty(0, T; L^2(\R^3)), \nabla\psi_t\in L^2(Q_T),\psi(0)=\u_0\}.
\end{split}
\end{equation*}
Consider the following linearized problem:
\begin{subequations}\label{e3}
\begin{align}
&\partial_t\r+\Dv(\r v)=0, \label{e31}\\
&\r\partial_t\u-\mu\D\u-(\mu+\lambda)\nabla\Dv\u=-\r v\cdot\nabla v-\nabla P+\Dv(\r \F\F^\top),  \label{e32}\\
&\partial_t \F+v\cdot\nabla\F=\nabla v \,\F. \label{e33}
\end{align}
\end{subequations}
with the given $v\in\mathcal{W}$  and the initial condition \eqref{IC}.

\subsection{Solvability of the density with a fixed velocity}
Let $A_j(x, t)$, $j=1, ..., n$,  be symmetric $m\times m$
matrices, $B(x,t)$ an $m\times m$ matrix, $f(x,t)$ and $V_0(x)$
two $m$-dimensional vector functions defined in $\R^n\times(0,T)$
and $\R^n$, respectively.

For the Cauchy problem of the linear system on $V\in\R^m$:
\begin{equation}\label{1000}
\begin{cases}
&\displaystyle \partial_tV+\sum_{j=1}^nA_j(x,t)\partial_{x_j}V+B(x,t)V=f(x,t),\\
&V(x,0)=V_0(x),
\end{cases}
\end{equation}
we have

\begin{Lemma}\label{l1}
Assume that $$A_j\in [C(0,T; H^s(\R^n))\cap C^1(0,T;
H^{s-1}(\R^{n}))]^{m\times m}, \;  j=1,...,n,$$
$$B\in C((0,T), H^{s-1}(\R^n))^{m\times m},\quad f\in C((0,T), H^s(\R^n))^m,\quad V_0\in H^s(\R^n)^m,$$
with $s>\f{n}{2}+1$  an integer. Then there exists a unique
solution to \eqref{1000}, i.e, a function
$$V\in [C([0,T), H^s(\R^n))\cap C^1((0,T), H^{s-1}(\R^n))]^m$$
satisfying \eqref{1000} pointwise (i.e. in the classical sense).
\end{Lemma}
\begin{proof}
This lemma is a direct consequence of Theorem 2.16 in \cite{AI}
with $A_0(x,t)=I$.
\end{proof}

To solve the density with respect to the velocity, we have
\begin{Lemma}\label{r}
Under the same conditions as Theorem \ref{T20}, there is a unique
strictly positive function
$$\r:=\mathcal{S}(v)\in W^{1,2}(0,T; L^q(\R^3)\cap L^2(\R^3))\cap L^\infty(0,T; W^{1,q}(\R^3)\cap H^1(\R^3))$$
which satisfy the continuity equation \eqref{e31}. Moreover, the density satisfies the following estimate:
$$\|\nabla \r\|^q_{L^\infty(0,T; L^q(\R^3))\cap L^2(\R^3)}\le \left(\|\nabla\r_0\|^q_{L^q(\R^3)\cap L^2(\R^3)}
+\sqrt{T}\|v\|_{\mathcal{W}}\right)\textrm{exp}\left(C\sqrt{T}\|v\|_{\mathcal{W}}\right).$$
\end{Lemma}

Here, and in what follows, the notation $C$ stands for a generic
positive constant, and in some cases, we will specify its
dependence on parameters by the notation $C(\cdot)$; and
$W^{1,2}(0,T;X)=\{f:f,f_t\in L^2(0,T,X)\}$. In many cases below, we drop the constant $C$ for the simplicity of notations.

\begin{proof}
The proof of the first part of this lemma is similar to that of
Lemma \ref{g} below, and can also be found in Theorem 9.3 in
\cite{AI}. The positivity of the density follows directly from the
observation: by writing \eqref{e31} along characteristics
$\f{d}{dt}X(t)=v,$
$$\f{d}{dt}\r(t,X(t))=-\r(t, X(t))\Dv v(t,X(t)), \quad X(0)=x,$$
and with the help of Gronwall's inequality,
\begin{equation*}
\begin{split}
\alpha\textrm{exp}(-\sqrt{T}\|v\|_{\mathcal{W}})&\le(\inf_{x}\r_0)\textrm{exp}\left(-\int_0^t\|\Dv
v(t)\|_{L^\infty(\R^3)}ds\right)\le \r(t,x)\\&\le(\sup_{x}\r_0)
\textrm{exp}\left(\int_0^t\|\Dv
v(t)\|_{L^\infty(\R^3)}ds\right)\le\beta\textrm{exp}(\sqrt{T}\|v\|_{\mathcal{W}}).
\end{split}
\end{equation*}

Now, we can assume that the continuity equation holds pointwise
in the following form
$$\partial_t \r+\r\, \Dv v+v\cdot \nabla \r=0.$$
Taking the gradient in both sides of the above identity,
multiplying by $|\nabla \r|^{q-2}\nabla \r$ and then integrating
over $\R^3$, we get, by Young's inequality,
\begin{equation}\label{y11}
\begin{split}
\f{1}{q}\f{d}{dt}\|\nabla \r\|^q_{L^q(\R^3)}&\le
\int_{\R^3}|\nabla \r|^q|\Dv v|dx+\int_{\R^3} \r|\nabla \r|^{q-1}
|\nabla\Dv v|dx\\&\qquad+\int_{\R^3}|\nabla v||\nabla
\r|^qdx-\f{1}{q}
\int_{\R^3} v\nabla|\nabla \r|^qdx\\
&\le\|\nabla \r\|^q_{L^q}\left(\|\nabla v\|_{L^\infty}+\|\r\|_{L^\infty}\|\nabla\Dv v\|_{L^q}\right)\\
&\quad+\f{1}{q}\int_{\R^3}\Dv v|\nabla \r|^qdx
+\|\r\|_{L^\infty}\|\nabla\Dv v\|_{L^q}\\
&\le C\|\nabla
\r\|^q_{L^q}\|v\|_{W^{2,q}}+\|\r\|_{L^\infty}\|\nabla\Dv
v\|_{L^q},
\end{split}
\end{equation}
since $W^{1,q}(\R^3)\hookrightarrow L^\infty(\R^3)$ as $q>3$.
Using Gronwall's inequality, we conclude that
\begin{equation}\label{y12}
\begin{split}
\|\nabla \r(t)\|^q_{L^q(\R^3)}&\le
\left(\|\nabla\r_0\|^q_{L^q}+\int_0^t\|\r\|_{L^\infty}\|\nabla\Dv
v\|_{L^q}ds\right)\textrm{exp}\left(\int_0^t\|v\|_{W^{2,q}}ds\right)\\
&\le\left(\|\nabla\r_0\|^q_{L^q}+\sqrt{t}\|v\|_{\mathcal{W}}\right)\textrm{exp}\left(\sqrt{t}\|v\|_{\mathcal{W}}\right).
\end{split}
\end{equation}
The proof is complete.
\end{proof}

\subsection{Solvability of the deformation gradient with a fixed velocity}

Due to the hyperbolic structure of \eqref{e33}, we can apply Lemma
\ref{l1} again to solve the deformation gradient $\F$ in terms of
the given velocity. For this purpose, we have
\begin{Lemma}\label{g}
Under the same conditions as Theorem \ref{T20}, there is a unique
function
$$\F:=\mathcal{T}(v)\in W^{1,2}(0,T; L^q(\R^3)\cap L^2(\R^3))\cap L^\infty(0,T; W^{1,q}(\R^3)\cap W^{1,2}(\R^3))$$
which satisfies the equation \eqref{e33}. Moreover, the
deformation gradient satisfies
$$\|\F\|_{L^\infty(0,T; W^{1,q}(\R^3))\cap H^1(\R^3)}\le \left(\|
\F(0)\|_{W^{1,q}(\R^3)\cap H^1(\R^3)
}+\sqrt{T}\|v\|_{\mathcal{W}}\right)\textrm{exp}(\sqrt{T}\|v\|_{\mathcal{W}}).$$
\end{Lemma}
\begin{proof}
First, we assume that $v\in C^1(0,T; C_0^\infty(\R^3)), \quad
\F_0\in C_0^\infty(\R^3))$. Then, we can rewrite \eqref{e33} in
the component form of columns as
$$\partial_t \F_k+v\cdot\nabla \F_k=\nabla v \F_k,\quad \textrm{for all}\quad 1\le k\le 3.$$
Applying Lemma \ref{l1} with $A_j(x,t)=v_j(x,t) I$ for $1\le j\le
3$, $B(x,t)=-\nabla v$, and $f(x,t)=0$, we get a solution
$$\F\in \bigcap_{s=3}^\infty\big\{C^1(0,T, H^{s-1}(\R^3))\cap C(0,T; H^s(\R^3))\big\}.$$
This implies, by the Sobolev imbedding theorems,
$$\F\in \bigcap_{k=1}^\infty C^1(0,T; C^k(\R^3))=C^1(0,T; C^\infty(\R^3)).$$

Next, for $v\in \mathcal{W}$, by an argument of dense sets, there is a sequence of
functions $v_n$ in the space $C^1(0,T; C_0^\infty(\R^3))$, $v_n\rightarrow
v$ in $\mathcal{W}$,  and $\F^n_0\in C_0^\infty(\R^3)$,
$\F^n_0\rightarrow \F_0$ in $W^{1,q}(\R^3)\cap H^1(\R^3)$. Hence,
$v_n\rightarrow v$ in $C(B(0,a)\times(0,T))$ for all $a>0$ where
$B(0,a)$ denotes the ball with radius $a$ and centered at the
origin. According to the previous result, there are
$\{\F_n\}_{n=1}^\infty$ satisfying
\begin{equation}\label{2}
\partial_t \F_n+v_n\cdot\nabla \F_n=\nabla v_n\F_n,
\end{equation}
with $\F_n(0)=\F^n_0$, $\F_n\in C^1(0,T; C^\infty(\R^3))$.
Multiply \eqref{2} by $|\F_n|^{p-2}\F_n$ for any $p\ge 2$, and
integrating over $\R^3$, by integration by parts, we obtain, using
Young's inequality,
\begin{equation*}
\begin{split}
\f{1}{p}\f{d}{dt}\int_{\R^3}|\F_n|^p dx&=-\f{1}{p}\int_{\R^3}
v_n\cdot\nabla |\F_n|^p dx
+\int_{\R^3}\nabla v_n|\F_n|^{p-2}\F^2_n dx\\
&\le \f{1+p}{p}\|\F_n\|_{L^p}^p\|\nabla v_n\|_{L^\infty}.
\end{split}
\end{equation*}
Then, by Gronwall's inequality, one has
\begin{equation*}
\begin{split}
\int_{\R^3}|\F_n|^p dx(t)& \le\int_{\R^3}|\F_n(0)|^p
dx\textrm{exp}\left(\int_0^t(p+1)\|\nabla
v_n\|_{L^\infty}ds\right)\\
&\le\int_{\R^3}|\F_n(0)|^p dx\textrm{exp}\left(\int_0^t(p+1)\|
v_n\|_{W^{2,p}}ds\right).
\end{split}
\end{equation*}
Thus,
\begin{equation*}
\begin{split}
\|\F_n\|_{L^\infty(0,T; L^p(\R^3))}&\le
\textrm{exp}\left(\f{p+1}{p}\|v_n\|_{\mathcal{W}}\sqrt{t}\right)\|\F_n(0)\|_{L^p(\R^3)}\\
&\le\textrm{exp}\left(2\|v_n\|_{\mathcal{W}}\sqrt{t}\right)\|\F_n(0)\|_{L^p(\R^3)}
<\infty.
\end{split}
\end{equation*}
Letting $p\rightarrow\infty$, one obtains
$$\|\F_n\|_{L^\infty(Q_T)}\le
\textrm{exp}(\|v_n\|_{\mathcal{W}}\sqrt{t})\|\F_n(0)\|_{L^\infty(\R^3)}\le\textrm{exp}(\|v_n\|_{\mathcal{W}}\sqrt{t})\|\F_n(0)\|_{W^{1,q}(\R^3)}
<\infty.$$ Hence, up to a subsequence, we can assume that $v_n$
were chosen so that
$$\F_n\rightarrow \F\quad \textrm{weak-* in}\quad L^\infty(0,T; L^q(\R^3)).$$

Taking the gradient in both sides of \eqref{2}, multiplying by
$|\nabla \F_n|^{q-2}\nabla \F_n$ and then integrating over $\R^3$,
we get, with the help of H\"{o}lder's inequality and Young's inequality,
\begin{equation}\label{y3}
\begin{split}
&\f{1}{q}\f{d}{dt}\|\nabla \F_n\|^q_{L^q(\R^3)}\\
&\le 2\int_{\R^3}|\nabla \F_n|^q|\nabla v_n|dx +\int_{\R^3}
|\F_n||\nabla \F_n|^{q-1}|\nabla\nabla v_n|dx-\f{1}{q}
\int_{\R^3} v_n\nabla|\nabla \F_n|^qdx\\
&\le C\int_{\R^3}|\nabla \F_n|^q|\nabla v_n|dx+\int_{\R^3} |\F_n||\nabla \F_n|^{q-1}|\nabla\nabla v_n|dx\\
&\le C\|\nabla \F_n\|^q_{L^q}\|v_n\|_{W^{2,q}}+\|\F_n\|_{L^\infty}\|v_n\|_{W^{2,q}}\|\nabla \F_n\|_{L^q}^{q-1}\\
&\le C\|\nabla
\F_n\|^q_{L^q}\|v_n\|_{W^{2,q}}+C\|v_n\|_{W^{2,q}}\|\nabla\F_n\|_{L^q}^{q-1},
\end{split}
\end{equation}
since $q>3$. Using Gronwall's inequality, we conclude that
\begin{equation*}
\begin{split}
\|\nabla \F_n(t)\|_{L^q(\R^3)}&\le \left(\|\nabla
\F_n(0)\|_{L^q}+\int_0^t\|v_n\|_{W^{2,q}}ds\right)\textrm{exp}(\int_0^t\|v_n\|_{W^{2,q}}ds)\\
&\le \left(\|\nabla
\F_n(0)\|_{L^q}+\sqrt{t}\|v_n\|_{\mathcal{W}}\right)\textrm{exp}(\sqrt{t}\|v_n\|_{\mathcal{W}}),
\end{split}
\end{equation*}
and hence,
\begin{equation}\label{y31}
\begin{split}
\|\nabla \F\|_{L^\infty(0,T; L^q(\R^3))}&\le\liminf_{n\rightarrow\infty}\|\nabla \F_n\|_{L^\infty(0,T; L^q)}\\
&\le \left(\|\nabla
\F(0)\|_{L^q}+\sqrt{T}\|v\|_{\mathcal{W}}\right)\textrm{exp}(\sqrt{T}\|v\|_{\mathcal{W}}).
\end{split}
\end{equation}
Thus,
$$\|\F\|_{L^\infty(0,T; W^{1,q}(\R^3))}\le \left(\|
\F(0)\|_{W^{1,q}}+\sqrt{T}\|v\|_{\mathcal{W}}\right)\textrm{exp}(\sqrt{T}\|v\|_{\mathcal{W}})<\infty.$$

Passing to the limit as $n\rightarrow\infty$ in \eqref{2}, we show
that \eqref{e33} holds at least in the sense of distributions.
Therefore, $\partial_t \F\in L^2(0,T; L^2(\R^3))$, then $\F\in
W^{1,2}(0,T; L^q(\R^3)\cap L^2(\R^3))$.
The proof is complete.
\end{proof}

\subsection{Local solvability of \eqref{e32}}

For simplicity of the presentation, we consider the case $\mu=1$
and $\lambda=0$ without loss of generality. In order to solve
\eqref{e32}, we consider the bilinear form $E(\u,\psi)$ and linear
functional $L(\psi)$ defined by
\begin{equation*}
\begin{split}
E(\u,\psi)&=\int_0^T(\r\partial_t\u-\D\u-\nabla\Dv\u,
\partial_t\psi-k(\D\psi+\nabla\Dv\psi))dt\\
&\qquad-(\u(0), \D\psi(0)+\nabla\Dv\psi(0)),
\end{split}
\end{equation*}
\begin{equation*}
\begin{split}
L(\psi)&=-\int_0^T(\r v\cdot\nabla v+\nabla P-\Dv(\r\F\F^\top),
\partial_t\psi-k(\D\psi+\nabla\Dv\psi))dt\\&\qquad-(\u_0,
\D\psi(0)+\nabla\Dv\psi(0)),
\end{split}
\end{equation*}
with $$k=(2\|\r\|_{L^\infty(Q_T)})^{-1}$$ for $\psi\in
\mathcal{W'}$, where $(\cdot,\cdot)$ denotes the inner product in
$L^2$.

We first notice that $L(\psi)$ is a linear continuous functional
of $\psi$ with respect to the norm $\|\psi\|_{\mathcal{W}'}$.
Moreover, we have
\begin{equation*}
\begin{split}
&E(\psi,\psi)\\&=\int_0^T\left(\|\sqrt{\r}\partial_t\psi\|^2_{L^2}+k\|\D\psi+\nabla\Dv\psi\|^2_{L^2}-k(\r\partial_t\psi,
\D\psi+\nabla\Dv\psi\right)dt\\
&\quad+\f{1}{2}\left(\|\nabla\psi(T)\|^2_{L^2}+\|\nabla\psi(0)\|^2_{L^2}+\|\Dv\psi(T)\|^2_{L^2}+\|\Dv\psi(0)\|^2_{L^2}\right)\\
&\ge\int_0^T\left(\|\sqrt{\r}\partial_t\psi\|^2_{L^2}+k\|\D\psi+\nabla\Dv\psi\|^2_{L^2}-\f{3}{4}\|\sqrt{\r}\partial_t\psi\|^2_{L^2}-\f{2k}{3}\|\D\psi+\nabla\Dv\psi\|^2_{L^2}\right)dt\\
&\quad+\f{1}{2}\left(\|\nabla\psi(T)\|^2_{L^2}+\|\nabla\psi(0)\|^2_{L^2}+\|\Dv\psi(T)\|^2_{L^2}+\|\Dv\psi(0)\|^2_{L^2}\right)\\
&\ge c_0\|\psi\|_{\mathcal{W'}}^2,
\end{split}
\end{equation*}
for some $c_0>0$, since
$$\|\D\psi+\nabla\Dv\psi\|_{L^2}\ge c_0\|\psi\|_{H^2}$$ from the
theory of elliptic operators. Hence, by the Lax-Milgram theorem
(see \cite{GT}), there exists a $\u\in \mathcal{W'}$ such that
\begin{equation}\label{cl6}
E(\u,\psi)=L(\psi)
\end{equation}
for every $\psi\in\mathcal{W'}$.

Now, let $\ov{\psi}$ be the solution of the problem
\begin{equation*}
\begin{split}
&\partial_t \ov{\psi}-k(\D\ov{\psi}+\nabla\Dv\ov{\psi})=0,\\
&\ov{\psi}(0)=h(x),
\end{split}
\end{equation*}
with $h(x)$ smooth enough. Replacing in \eqref{cl6} $\psi$ by
$\ov{\psi}$, one obtain
$$(\u(0)-\u_0, \D h+\nabla\Dv h)=0,$$ which implies $\u(0)=\u_0$.
Next, let $\widetilde{\psi}$ be a solution of the following
problem
\begin{equation*}
\begin{split}
&\partial_t \widetilde{\psi}-k(\D\widetilde{\psi}+\nabla\Dv\widetilde{\psi})=g(x,t),\\
&\widetilde{\psi}(0)=0,
\end{split}
\end{equation*}
with $g$ smooth enough. Replacing $\psi$ by $\widetilde{\psi}$ in
\eqref{cl6}, one obtain
$$\int_0^T(\r\partial_t\u-\D\u-\nabla\Dv\u+\r v\cdot\nabla
v-\Dv(\r \F\F^\top)+\nabla P, g)dt=0.$$ This implies that
$(\u,\r,\F)$ satisfies \eqref{e3} a.e. in $(0,T)\times\R^3$.

Next, we prove the higher regularity for $\u$; that is,
$\u\in\mathcal{W}$. First, we multiply \eqref{e32} by
$\partial_t\u$, and use integration by parts and Young's
inequality to obtain
\begin{equation}\label{cl7}
\begin{split}
\f{2}{3}&\int_0^t\|\sqrt{\r}\partial_t\u\|^2_{L^2}ds+\f{1}{2}\|\nabla\u(t)\|^2_{L^2}+\f{1}{2}\|\Dv\u(t)\|_{L^2}^2\\
&\le\f{1}{2}\|\nabla\u(0)\|^2_{L^2}+\f{1}{2}\|\Dv\u(0)\|_{L^2}^2-\int_0^t(\nabla
P,\partial_t\u)ds+\\&\quad+C\sup_{s\in(0,t)}
\|\r\|_{L^\infty}\int_0^t\|v\|_{L^6}^2\|\nabla
v\|_{L^3}^2ds+\int_0^t(\Dv(\r\F\F^\top),\partial_t \u)ds+\f{1}{4}\int_0^t\|\partial_t\nabla\u\|_{L^2}^2ds\\
&\le\f{1}{2}\|\nabla\u(0)\|^2_{L^2}+\f{1}{2}\|\Dv\u(0)\|_{L^2}^2+\int_0^t\|\nabla\r\|_{L^2}^2ds
+\f{1}{4}\int_0^t\|\partial_t\nabla\u\|_{L^2}^2ds+\\&\quad+C\sup_{s\in(0,t)}
\|\r\|_{L^\infty}\int_0^t\|v\|_{L^6}^2\|\nabla
v\|_{L^3}^2ds+\int_0^t\|\nabla\F\|_{L^2}^2ds+\f{1}{3}\int_0^t\|\sqrt{\r}\partial_t\u\|_{L^2}^2ds.
\end{split}
\end{equation}
In particular, if we multiply \eqref{e32} by $\partial_t\u$,
integrate over $\R^3$,  and let $t=0$, we obtain
\begin{equation*}
\begin{split}
\|\sqrt{\r}(0)\partial_t\u(0)\|_{L^2}^2
&=\int_{\R^3}\Big(\D\u(0)+\nabla\Dv\u(0)-\r(0)v(0)\cdot\nabla v(0)-\nabla P(\r(0))\\
&\qquad\qquad +\Dv(\r(0)\F(0)\F^\top(0))\Big)\partial_t\u(0) dx\\
&\le \|\sqrt{\r(0)}\|_{L^\infty}^{-1}\Big(\|\D\u_0\|_{L^2}+\|\r(0)\|_{L^\infty}\|v(0)\|_{L^\infty}\|\nabla v(0)\|_{L^2} \\
&\qquad\qquad\qquad\qquad +C\|\nabla\r(0)\|_{L^2}+\|\nabla\F(0)\|_{L^2}\Big)\|\sqrt{\r}(0)\partial_t\u(0)\|_{L^2},
\end{split}
\end{equation*}
which implies that
\begin{equation}\label{cl13}
\|\sqrt{\r}(0)\partial_t\u(0)\|_{L^2}\le
C(\|\D\u_0\|_{L^2}+\|v_0\|_{L^\infty}\|\nabla v_0\|_{L^2}
+C\|\nabla\r_0\|_{L^2} +\|\nabla\F_0\|_{L^2}).
\end{equation}

Now we differentiate \eqref{e32} with respect to $t$ so that we get
\begin{equation}\label{cl8}
\begin{split}
&\partial_t\r\partial_t\u+\r\partial^2_{tt}\u-\D\partial_t\u-\nabla\Dv(\partial_t\u)\\
&=-\partial_t\r v\cdot\nabla v-\r\partial_t v\cdot\nabla v-\r
v\cdot\nabla\partial_t v -\nabla\partial_t P+\partial_t\Dv(\r \F\F^\top).
\end{split}
\end{equation}
Multiplying \eqref{cl8} by $\partial_t\u$, integrating over
$\R^3$, and using the continuity equation, we obtain
\begin{equation}\label{cl9}
\begin{split}
&\f{1}{2}\f{d}{dt}\|\sqrt{\r}\partial_t\u\|^2_{L^2}+\f{1}{2}\int_{\R^3}\partial_t\r|\partial_t\u|^2dx+\|\nabla\partial_t\u\|^2_{L^2}+
\|\Dv\partial_t\u\|^2_{L^2}\\
&\le\|\r\|_{L^\infty}\|v\|_{L^\infty}\|\nabla v\|_{L^2}\|\nabla
v\|_{L^3}\|\partial_t\u\|_{L^6}+\|\r\|_{L^\infty}\|v\|_{L^6}^2\|\nabla(\nabla v)\|_{L^2}\|\partial_t\u\|_{L^6}\\
&\quad+\|\sqrt{\r}\|_{L^\infty}\|\sqrt{\r}\partial_t
v\|_{L^2}\|\nabla
v\|_{L^3}\|\partial_t\u\|_{L^6}+\|\r\|_{L^\infty}\|v\|_{L^\infty}\|v\|_{L^6}\|\nabla
v\|_{L^3}\|\nabla\partial_t\u\|_{L^2}\\
&\quad+\|\sqrt{\r}\|_{L^\infty}\|\nabla\partial_t
v\|_{L^2}\|v\|_{L^\infty}\|\sqrt{\r}\partial_t\u\|_{L^2}+\|\partial_t
P\|_{L^2}\|\nabla\partial_t\u\|_{L^2}\\
&\quad +\|\nabla\partial_t\u\|_{L^2}\Big(\|\nabla\r\|_{L^2}\|v\|_{L^\infty}\|\F\|_{L^\infty}^2+\|\r\|_{L^\infty}\|\nabla
v\|_{L^2}\|\F\|_{L^\infty}^2\\
&\qquad\qquad\qquad\qquad\qquad \qquad\qquad\qquad\qquad
+\|\r\|_{L^\infty}\|\F\|_{L^\infty}\|v\|_{L^\infty}\|\nabla\F\|_{L^2}\Big).
\end{split}
\end{equation}

Integrating \eqref{cl9} with respect to $t$, using the continuity
equation and the Gagliardo-Nirenberg inequality
$$\|\nabla\u\|_{L^3}^2\le C\|\nabla\u\|_{L^2}\|\D\u\|_{L^2},$$
we find, since $\r\in L^\infty((0,T)\times\R^3)$,
\begin{equation}\label{cl10}
\begin{split}
&\f{1}{2}\|\sqrt{\r}\partial_t\u(t)\|^2_{L^2}+\int_0^t(\|\nabla\partial_t\u\|^2_{L^2}+
\|\Dv\partial_t\u\|^2_{L^2})ds \\
&\le\f{1}{2}\|\sqrt{\r}\partial_t\u(0)\|^2_{L^2} \\
&\quad+\int_0^t\|\r\|_{L^\infty}\Big(\|v\|_{L^\infty}\|\nabla
v\|_{L^2}^{\f{3}{2}}\|\D
v\|_{L^2}^{\f{1}{2}}\|\nabla\partial_t\u\|_{L^2}+\|\nabla v\|_{L^2}^2\|\D v\|_{L^2}\|\nabla \partial_t\u\|_{L^2}\Big)ds\\
&\quad+\int_0^t\|\sqrt{\r}\|_{L^\infty}\Big(\|\sqrt{\r}\partial_t
v\|_{L^2}\|\nabla
v\|_{L^2}^{\f{1}{2}}\|\D \u\|_{L^2}^{\f{1}{2}}\|\nabla\partial_t\u\|_{L^2}+\|v\|_{L^\infty}\|\nabla\partial_t\u\|_{L^2}\|\sqrt{\r}\partial_t\u\|_{L^2}\Big)ds\\
&\quad+\int_0^t\|\sqrt{\r}\|_{L^\infty}\|\nabla\partial_t
v\|_{L^2}\|v\|_{L^\infty}\|\sqrt{\r}\partial_t\u\|_{L^2}ds+\int_0^t\|\partial_t
P\|_{L^2}\|\nabla\partial_t\u\|_{L^2}ds\\&\quad+C_1\int_0^t\|\nabla\partial_t\u\|_{L^2}(\|v\|_{L^\infty}+\|\nabla v\|_{L^2})ds\\
&\le\f{1}{2}\|\sqrt{\r}\partial_t\u(0)\|^2_{L^2}+C_{\dl}\int_0^t\|v\|_{L^\infty}^2\|\sqrt{\r}\partial_t\u\|_{L^2}^2ds+
C_{\dl}\int_0^t(\|v\|_{L^\infty}^2\|\nabla v\|_{L^2}^3\|\D v
\|_{L^2}\\&\quad+C_{\dl}\|\nabla v\|^4_{L^2}\|\D
v\|_{L^2}^2)ds+C_{\dl}\int_0^t\|\sqrt{\r}\partial_t
v\|_{L^2}^2\|\nabla v\|_{L^2}\|\D
v\|_{L^2}ds+C_{\dl}\int_0^t\|\Dv(\r v)
\|_{L^2}^2ds\\&\quad+\dl\int_0^t(\|\nabla\partial_t\u\|_{L^2}^2+\|\sqrt{\r}\partial_t\u\|^2_{L^2})ds+C_{\dl}\int_0^t(\|v\|_{L^\infty}^2+\|\nabla
v\|_{L^2}^2)ds+\int_0^t\|\nabla\partial_t v\|_{L^2}^2ds,
\end{split}
\end{equation}
where $\dl$ is a small constant, $C_1$ depends on $v$,  $\r_0$,
and $\F_0$, since
 $$\nabla \r\in L^\infty(0,T; L^q(\R^3)\cap L^2(\R^3)),\quad \nabla \F\in L^\infty(0,T; L^q(\R^3)\cap L^2(\R^3))$$
  as stated in Lemma \ref{r} and Lemma \ref{g}. Since we
are only interested in the local existence, so we can restrict
$t\le\ov{T}\le 1$.

Adding \eqref{cl7} and \eqref{cl10} for some suitable
$\dl<\f{1}{2}$, one obtains, first by Gronwall's inequality,
$$\|\sqrt{\r}\partial_t\u\|_{L^\infty(0,t; L^2(\R^3))}\le C,$$
and, second,
\begin{equation*}
\begin{split}
&\f{1}{2}\|\sqrt{\r}\partial_t\u(t)\|^2_{L^2}+\f{1}{4}\int_0^t(\|\nabla\partial_t\u\|^2_{L^2}+
\|\Dv\partial_t\u\|^2_{L^2})ds+\f{1}{3}\int_0^t\|\sqrt{\r}\partial_t\u\|^2_{L^2}ds\\&\quad+\f{1}{2}\|\nabla\u(t)\|^2_{L^2}+\f{1}{2}\|\Dv\u(t)\|_{L^2}^2\le
C,
\end{split}
\end{equation*}
which implies
\begin{equation}\label{cl11}
\begin{cases}
&\sqrt{\r}\partial_t\u\in L^\infty(0,T; L^2(\R^3));\quad
\partial_t\u\in L^2(0,T; H^1(\R^3));\\
&\nabla\u\in L^\infty(0,T; L^2(\R^3)).
\end{cases}
\end{equation}

On the other hand, we can rewrite \eqref{e32} as
$$-\D\u-\nabla\Dv\u=-\r\partial_t\u-\r v\cdot\nabla v-\nabla
P+\Dv(\r\F\F^\top),$$ which is a strongly elliptic equation.
Hence, by the classical theory in the elliptic system, we deduce
that
$$\u\in L^\infty(0,T; H^2(\R^3)),$$
since from Lemmas \ref{r}-\ref{g} and \eqref{cl11} one has
$$-\r\partial_t\u-\r v\cdot\nabla v-\nabla P+\Dv(\r\F\F^\top)\in
L^\infty(0,T; L^2(\R^3)).$$ Moreover, since $\r$ is bounded from
below and $\sqrt{\r}\partial_t\u\in L^2(0,T; L^2(\R^3))$, we know
that $\partial_t\u\in L^2(0,T; L^2(\R^3))$. Hence, by
the Gagliardo-Nirenberg inequality, as $q\in (3,6]$,
$$\|\partial_t\u\|_{L^2(0,T; L^q(\R^3))}\le
C\|\partial_t\u\|_{L^2((0,T)\times\R^3)}^\theta\|\nabla\partial_t\u\|_{L^2((0,T)\times\R^3)}^{1-\theta},$$
for some $\theta\in [0,1)$. This implies that, by \eqref{cl11},
$\partial_t\u\in L^2(0,T; L^q(\R^3))$. Thus, by the classical
elliptic theory, we obtain $\u\in L^2(0,T; W^{2,q}(\R^3))$ since
now $$-\r\partial_t\u-\r v\cdot\nabla v-\nabla
P+\Dv(\r\F\F^\top)\in L^2(0,T; L^q(\R^3)).$$
Hence, we can conclude that $\u\in\mathcal{W}$.

\subsection{Existence for \eqref{e1}}
The above argument  leads us to define the map
$$\u=\mathcal{H}(v)$$
from $\mathcal{W}$ to itself through the maps
$g:v\mapsto\mathcal{S}(v)$, $f:v\mapsto\mathcal{T}(v)$ and $d:
(\mathcal{S}(v), v,\mathcal{T}(v))\mapsto\u$. Hence the solution
of \eqref{e1} is obtained from a fixed point of the map $\mathcal{H}$.
To find a fixed point of  $\mathcal{H}$,  we will use the Schauder-Tychonoff fixed point theorem
(Theorem 5.28, \cite{RW}). Define
\begin{equation*}
\begin{split}
M=\Big\{\psi: & \; \max\big(\|\psi\|_{L^2(0,T; W^{2,q}(\R^3)\cap H^2(\R^3))}, \,
 \|\sqrt{\r}\partial_t\psi\|_{L^\infty(0,T;L^2(\R^3))}, \\&\qquad\qquad \|\psi\|_{L^\infty(0,T; H^2(\R^3))}, \,
\|\partial_t\psi\|_{L^2(0,T; H^1(\R^3))}\big)\le r\Big\},
\end{split}
\end{equation*}
where $$r=C_0(\|\u_0\|_{H^2}+\|\u_0\|_{L^\infty}\|\nabla
\u_0\|_{L^2} +\|\r_0\|_{W^{1,q}\cap H^1} +\|\F_0\|_{W^{1,q}\cap
H^1})$$ with some sufficiently large $C_0>0$. Clearly, $M$ is a
compact and convex set in $L^2((0,T)\times\R^3)$. Hence, we need
to show that $\mathcal{H}(M)\subseteq M$ (i.e., $\mathcal{H}$ maps
$M$ into $M$) and $\mathcal{H}$ is continuous in $M$ with respect
to the norm in $L^2((0,T)\times\R^3)$.

We first prove that $\mathcal{H}(M)\subset M$ for some $T=\ov{T}$.
Indeed, assuming $v\in M$, from Lemma \ref{r} and Lemma \ref{g},
we know that
\begin{equation}\label{cl12}
\begin{cases}
\alpha\textrm{exp}(-r\sqrt{t})\le\r(x,t)\le \beta
\textrm{exp}(r\sqrt{t});\\
\|\r\|_{L^\infty(0,t; W^{1,q}(\R^3))}\le
\left(\|\r_0\|_{W^{1,q}(\R^3)}
+\sqrt{t}r\right)\textrm{exp}\left(\sqrt{t}r\right);\\
\|\F\|_{L^\infty(0,t; W^{1,q}(\R^3))}\le \left(\|
\F(0)\|_{W^{1,q}(\R^3)}+\sqrt{t}r\right)\textrm{exp}(\sqrt{t}r).
\end{cases}
\end{equation}
Hence, from \eqref{cl10} and \eqref{cl12}, it follows that
\begin{equation*}
\begin{split}
&\f{1}{2}\|\sqrt{\r}\partial_t\u(t)\|_{L^2}^2+\f{1}{2}\int_0^t\|\nabla\partial_t\u\|_{L^2}^2ds \\
&\le \f{1}{2}\|\sqrt{\r}\partial_t\u(0)\|_{L^2}^2 +C_{\dl}t(r^2+\dl)\|\sqrt{\r}\partial_t\u\|_{L^\infty(0,t;
L^2(\R^3))}^2+C_{\dl}t(r^4+r^6)+Ctr^2.
\end{split}
\end{equation*}
Using \eqref{cl13}, and taking $\dl$ and $\ov{T}$ sufficiently
small, we derive from the above inequality that
\begin{equation*}
\|\sqrt{\r}\partial_t\u\|_{L^\infty(0,\ov{T};
L^2(\R^3))}^2+\|\nabla\partial_t\u\|_{L^2(0,\ov{T};
L^2(\R^3))}^2\le \f{1}{3}r^2.
\end{equation*}
Since $\r$ is bounded from below, we obtain
\begin{equation}\label{cl14}
\|\sqrt{\r}\partial_t\u\|_{L^\infty(0,\ov{T};
L^2(\R^3))}^2+\|\partial_t\u\|_{L^2(0,\ov{T}; H^1(\R^3))}^2\le
\f{2}{3}r^2.
\end{equation}

To estimate the norm $\|\u\|_{L^\infty(0,T; H^2(\R^3))}$, noticing
that
$$-\D\u-\nabla\Dv\u=-\r\partial_t\u-\r v\cdot\nabla v-\nabla
P(\r)+\Dv(\r \F\F^\top),$$ one has from the theory for elliptic
equations that
\begin{equation*}
\begin{split}
\|\u\|_{L^\infty(0,\ov{T};H^2(\R^3))}&\le
\|\sqrt{\r}\|_{L^\infty}\|\sqrt{\r}\partial_t\u\|_{L^\infty(0,\ov{T};
L^2)}+\|\r\|_{L^\infty}\|v\|_{L^\infty}\|\nabla
v\|_{L^\infty(0,\ov{T};L^2)}\\&\quad+C\|\nabla\r\|_{L^\infty(0,\ov{T};
L^2)}+\|\nabla\F\|_{L^\infty(0,\ov{T}; L^2)};
\end{split}
\end{equation*}
and hence
$$\|\D\u+\nabla\Dv\u\|^2_{L^\infty(0,\ov{T};H^2)}\le r^2.$$

Next, we need  to obtain estimates on $\|\u\|_{L^2(0,T;
W^{2,q}(\R^3))}$. Indeed, we have, by the classical theory of
elliptic equations,
\begin{equation*}
\begin{split}
\int_0^{\ov{T}}\|\u\|_{W^{2,q}}^2dt&\le
C\int_0^{\ov{T}}(\|v\|^2_{L^\infty}\|\nabla
v\|_{L^q}^2+\|\nabla\r\|_{L^q}^2+\|\partial_t\u\|_{L^q}^2+\|\nabla\F\|_{L^q}^2)dt\\
&\le Cr^4\ov{T}+r^2\ov{T}\le r^2,
\end{split}
\end{equation*}
for some sufficiently small $\ov{T}$. Hence, we show that
$\mathcal{H}(M)\subseteq M$.

Finally, we need to prove the continuity of $\mathcal{H}$ in $M$.
First we observe that if $\{v_n\}_{n=1}^\infty\subseteq M$, then
there exists a subsequence (still denoted by
$\{v_n\}_{n=1}^\infty$) such that as $n\rightarrow \infty$,
$v_n\rightarrow v$ strongly in $M$. Let $\r_n$ and $\r$ be the
solutions of
$$\partial_t\r_n+\Dv(\r_n v_n)=0,$$
and
$$\partial_t\r+\Dv(\r v)=0,$$
with $\r_n(0)=\r(0)=\r_0$, respectively. Denoting
$\ov{\r_n}=\r_n-\r$, then $\ov{\r_n}$ satisfies
$$\partial_t\ov{\r_n}+v_n\cdot\nabla\ov{\r_n}+(v_n-v)\cdot\nabla\r+\ov{\r_n}\Dv
v_n+\r\Dv(v_n-v)=0, $$ with $\ov{\r_n}(0)=0$. Repeating the
argument in Lemma \ref{r}, we have
$$\|\ov{\r_n}\|_{L^2}^2\le
\textrm{exp}(Cr\ov{T})\int_0^{\ov{T}}\left(\|(v-v_n)\cdot\nabla\r\|_{L^2}^2+\|\r\Dv(v_n-v)\|_{L^2}^2\right)dt,$$
which implies that $\r_n\rightarrow\r$ strongly in
$L^\infty(0,\ov{T}; L^2(\R^3))$. Similarly, we can show that
$\F_n\rightarrow\F$ strongly in $L^\infty(0,\ov{T}; L^2(\R^3))$.
Now let $\u_n$ and $\u$ be the solutions of \eqref{e32}
corresponding to $v_n$ and $v$ with $\u_n(0)=\u(0)=\u_0$
respectively. Then one has,  denoting $U_n=\u_n-\u$ and
$V_n=v_n-v$,
\begin{equation*}
\begin{split}
&\r_n\partial_t U_n-\D U_n-\nabla\Dv U_n \\
&=-\ov{\r_n}\partial_t\u-\r_n V_n\cdot\nabla v_n-\ov{\r_n}v\cdot\nabla v_n-\r v\cdot\nabla V_n\\
&\quad+\Dv(\r_n\F_n\F_n^\top-\r\F\F^\top)-\nabla P(\r_n)+\nabla
P(\r).
\end{split}
\end{equation*}
Multiplying the above equation by $\partial_t U_n$, integrating
over $(0,\ov{T})\times\R^3$, and thanking to the convergence of
$\r_n$ and $\F_n$, we can prove as a routine matter that $\nabla
U_n\rightarrow 0$ strongly in $L^2((0,\ov{T})\times\R^3)$ and
$\sqrt{\r_n}\partial_t U_n\rightarrow 0$ strongly in
$L^2((0,\ov{T})\times\R^3)$. Due to the lower bound of $\r_n$, we
deduce that $\partial_t U_n\rightarrow 0$ strongly in
$L^2((0,\ov{T})\times\R^3)$, and hence $U_n\rightarrow 0$ in
$L^2(0,\ov{T})\times\R^3)$ by using the identiy
$U_n(t)=\int_0^t\partial_t U_nds$ since $U_n(0)=0$. Thus, the map
$\mathcal{H}$ is continuous in $M$. The existence of a local
solution is completely proved.

\bigskip\bigskip

\section{Uniqueness}
In this section, we will prove the uniqueness of the solution
obtained in the previous section. Notice that, the argument in
Section 3 yields that $$\partial_t \u\in L^2(0,T; L^2\cap
L^q(\R^3), \; \nabla \r\in L^2(0,T; L^2\cap L^q(\R^3)), \; \nabla
\F\in L^2(0,T; L^2\cap L^q(\R^3)),$$ for $q>3$.  Hence, using the
interpolation, we see that
$$\partial_t\u\in L^{2}(0,T; L^3(\R^3),\; \nabla \r\in L^2(0,T; L^3(\R^3)), \; \nabla \F\in L^2(0,T; L^3(\R^3)).$$
Now, assume that $\u_1$, $\u_2$ satisfy \eqref{e1} for some
$T>0$, and let 
$$r:=\mathcal{S}(\u_1)-\mathcal{S}(\u_2),\quad
v:=\u_1-\u_2,\quad G:=\mathcal{T}(\u_1)-\mathcal{T}(\u_2).$$
 Then, we have
\begin{equation}\label{12}
\begin{cases}
&\partial_t r+\u_1\cdot\nabla r+v\cdot\nabla \mathcal{S}(\u_2)+r\Dv \u_1+\mathcal{S}(\u_2)\Dv v=0,\\
&r(0)=0.
\end{cases}
\end{equation}
Multiplying \eqref{12} by $r$, and integrating over $\R^3$, we get
\begin{equation*}
\begin{split}
&\f{1}{2}\f{d}{dt}\|r\|_{L^2}^2-\f{1}{2}\int_{\R^3}|r|^2\Dv
\u_1dx+\int_{\R^3} v\nabla \mathcal{S}(\u_2)r dx
\\&\quad+\int_{\R^3}|r|^2\Dv \u_1dx+\int_{\R^3} r \mathcal{S}(\u_2)\Dv
vdx =0,
\end{split}
\end{equation*}
which yields
\begin{equation}\label{1213}
\begin{split}
\f{d}{dt}\|r\|^2_{L^2(\R^3)}&\le \|\Dv
\u_1\|_{L^\infty}\|r\|^2_{L^2}+\varepsilon\|\nabla v\|_{L^2}^2
+C(\varepsilon)\|\nabla
\mathcal{S}(\u_2)r\|^2_{L^{\f{6}{5}}}\\&\qquad+\varepsilon\|\nabla
v\|^2_{L^2}
+C(\varepsilon)\|\mathcal{S}(\u_2)\|_{L^\infty}^2\|r\|^2_2\\
&\le \|\Dv \u_1\|_{L^\infty}\|r\|^2_{L^2}+\varepsilon\|\nabla
v\|_{L^2}^2+C(\varepsilon) \|\nabla
\mathcal{S}(\u_2)\|^2_{L^3}\|r\|_{L^2}^2\\&\qquad+\varepsilon\|\nabla
v\|^2_{L^2}
+C(\varepsilon)\|\mathcal{S}(\u_2)\|_{L^\infty}^2\|r\|^2_2\\
&\le\eta_1(\varepsilon)\|r\|_{L^2}^2+2\varepsilon\|\nabla
v\|^2_{L^2},
\end{split}
\end{equation}
where $\eta_1(\varepsilon)=\|\Dv
\u_1\|_{L^\infty}+C(\varepsilon)(\|\nabla
\mathcal{S}(\u_2)\|^2_{L^3}+\|\mathcal{S}(\u_2)\|^2_{L^\infty}).$

Similarly, from \eqref{e33}, we obtain
\begin{equation}\label{1211}
\begin{cases}
&\partial_t G+\u_1\cdot\nabla G+v\cdot\nabla \mathcal{T}(\u_2)=\nabla \u_1 G+\nabla v \mathcal{T}(\u_2),\\
&G(0)=0.
\end{cases}
\end{equation}
Multiplying \eqref{1211} by $G$, and integrating over $\R^3$, we
get
\begin{equation*}
\begin{split}
&\f{1}{2}\f{d}{dt}\|G\|_{L^2}^2-\f{1}{2}\int_{\R^3}|G|^2\Dv \u_1dx+\int_{\R^3} v\cdot\nabla \mathcal{T}(\u_2):G dx\\
&=\int_{\R^3}G^\top\nabla \u_1 Gdx+\int_{\R^3} \nabla
v\mathcal{T}(\u_2): G dx,
\end{split}
\end{equation*}
which yields
\begin{equation}\label{1212}
\begin{split}
\f{d}{dt}\|G\|^2_{L^2(\R^3)}&\le \|\Dv
\u_1\|_{L^\infty}\|G\|^2_{L^2}+\varepsilon\|\nabla v\|_{L^2}^2
+C(\varepsilon)\|\nabla
\mathcal{T}(\u_2)G\|^2_{L^{\f{6}{5}}}\\&\qquad+\varepsilon\|\nabla
v\|^2_{L^2}
+C(\varepsilon)\|\mathcal{T}(\u_2)\|_{L^\infty}^2\|G\|^2_2\\
&\le \|\Dv \u_1\|_{L^\infty}\|G\|^2_{L^2}+\varepsilon\|\nabla
v\|_{L^2}^2+C(\varepsilon) \|\nabla
\mathcal{T}(\u_2)\|^2_{L^3}\|G\|_{L^2}^2\\&\qquad+\varepsilon\|\nabla
v\|^2_{L^2}
+C(\varepsilon)\|\mathcal{T}(\u_2)\|_{L^\infty}^2\|G\|^2_2\\
&\le\eta_2(\varepsilon)\|G\|_{L^2}^2+2\varepsilon\|\nabla
v\|^2_{L^2},
\end{split}
\end{equation}
where $\eta_2(\varepsilon)=\|\Dv
\u_1\|_{L^\infty}+C(\varepsilon)(\|\nabla
\mathcal{T}(\u_2)\|^2_{L^3} +\|\mathcal{T}(\u_2)\|^2_{L^\infty}).$

For each $\u_j$ $j=1, 2$, we deduce from \eqref{e32} that
\begin{equation*}
\begin{cases}
&\mathcal{S}(\u_j)\partial_t \u_j-\mu\Delta
\u_j-(\mu+\lambda)\nabla\Dv \u_j\\&\quad=-\mathcal{S}(\u_j)
(\u_j\cdot\nabla)\u_j-\nabla P(\mathcal{S}(\u_j))+\Dv(\mathcal{S}(\u_j)\mathcal{T}(\u_j)\mathcal{T}(\u_j)^\top),\\
&\u_j(0)=\u_0.
\end{cases}
\end{equation*}
Subtracting these equations, we obtain,
\begin{equation}\label{13}
\begin{split}
&\mathcal{S}(\u_1)\partial_t \u_1-\mathcal{S}(v_2)\partial_t
\u_2-\mu\D v-(\mu+\lambda)\nabla\Dv v\\&=
-\mathcal{S}(\u_1)(\u_1\cdot\nabla)\u_1+\mathcal{S}(\u_2)(\u_2\cdot\nabla)\u_2-\nabla
P(\mathcal{S}(\u_1))+\nabla
P(\mathcal{S}(\u_2))\\&\qquad+\Dv(\mathcal{S}(\u_1)\mathcal{T}(\u_1)\mathcal{T}(\u_1)^\top)-\Dv(\mathcal{S}(\u_2)\mathcal{T}(\u_2)\mathcal{T}(\u_2)^\top).
\end{split}
\end{equation}

Since
\begin{equation*}
\begin{split}
&-\mathcal{S}(\u_1)(\u_1\cdot\nabla)\u_1+\mathcal{S}(\u_2)(\u_2\cdot\nabla)\u_2\\&=-\mathcal{S}(\u_1)(v\cdot\nabla)\u_1-
(\mathcal{S}(\u_1)-\mathcal{S}(\u_2))(\u_2\cdot\nabla)\u_1-\mathcal{S}(\u_2)(\u_2\cdot\nabla)v,
\end{split}
\end{equation*}
and
\begin{equation*}
\begin{split}
&\mathcal{S}(\u_1)\mathcal{T}(\u_1)\mathcal{T}(\u_1)^\top-\mathcal{S}(\u_2)\mathcal{T}(\u_2)\mathcal{T}(\u_2)^\top\\&=\mathcal{S}(\u_1)G
\mathcal{T}(\u_1))^\top +r
\mathcal{T}(\u_2)\mathcal{T}(\u_1)^\top+\mathcal{S}(\u_2)\mathcal{T}(\u_2)G^\top,
\end{split}
\end{equation*}
we can rewrite \eqref{13} as
\begin{equation}\label{14}
\begin{split}
&\mathcal{S}(\u_1)\partial_t v-\mu\D v-(\mu+\lambda)\nabla\Dv
v\\&=-r\partial_t \u_2-\mathcal{S}(\u_1)(v\cdot\nabla)
\u_1-r(\u_2\cdot\nabla)\u_1-\mathcal{S}(\u_2)(\u_2\cdot\nabla)v-\nabla P(\mathcal{S}(\u_1))+\nabla P(\mathcal{S}(\u_2))\\
&\qquad+\Dv(\mathcal{S}(\u_1)G \mathcal{T}(\u_1)^\top+r
\mathcal{T}(\u_2)\mathcal{T}(\u_1)^\top+\mathcal{S}(\u_2)\mathcal{T}(\u_2)G^\top).
\end{split}
\end{equation}

Multiplying \eqref{14} by $v$, using the continuity equation
\eqref{e11} and integrating over $\R^3$, we deduce that
\begin{equation}\label{15}
\begin{split}
&\f{1}{2}\f{d}{dt}\int_{\R^3} \mathcal{S}(\u_1)|v|^2 dx+\int_{\R^3}(\mu|\nabla v|^2+(\lambda+\mu)|\Dv v|^2)dx\\
&=\int_{\R^3}
\Big\{\f{1}{2}\mathcal{S}(\u_1)(\u_1\cdot\nabla)v\cdot
v-r\partial_t \u_2v-\mathcal{S}(\u_1)(v\cdot\nabla)
\u_1v-r(\u_2\cdot\nabla)\u_1v\\&\qquad-\mathcal{S}(\u_2)(\u_2\cdot\nabla)vv-\nabla
P(\mathcal{S}(\u_1))v+\nabla
P(\mathcal{S}(\u_2))v\\&\qquad-\big(\mathcal{S}(\u_1)G \mathcal{T}(\u_1)^\top+r \mathcal{T}(\u_2)\mathcal{T}(\u_1)^\top+\mathcal{S}(\u_2)\mathcal{T}(\u_2)G^\top\big)\nabla v\Big\}dx\\
&\le \varepsilon \|\nabla
v\|^2_{L^2}+C(\varepsilon)\|\mathcal{S}(\u_1)\|^2_{L^\infty}
\|\u_1\|^2_{L^\infty}\|v\|^2_{L^2}+\varepsilon\|\nabla v\|^2_{L^2}+C(\varepsilon)\|\partial_t \u_2\|_{L^3}^2\|r\|_{L^2}^2\\
&\qquad +\|\mathcal{S}(\u_1)\|_{L^\infty}\|\nabla
\u_1\|_{L^\infty}\|v\|^2_{L^2}
+2\|\u_2\|_{L^\infty}\|\nabla \u_1\|_{L^\infty}(\|r\|_{L^2}^2+\|v\|^2_{L^2})\\
&\qquad +\varepsilon\|\nabla
v\|^2_{L^2}+C(\varepsilon)\|\mathcal{S}(\u_2)\|_{L^\infty}^2
\|\u_2\|_{L^\infty}^2\|v\|_{L^2}^2+\varepsilon\|\nabla v\|_{L^2}^2\\
&\qquad+C(\varepsilon)(\sup\{P'(\eta): C(T)^{-1}\le\eta\le
C(T)\})^2\|r\|^2_{L^2}
+\varepsilon\|\nabla v\|^2_{L^2}\\
&\qquad+C(\varepsilon)(\|\mathcal{S}(\u_1)\|_{L^\infty}^2\|\mathcal{T}(\u_1)\|^2_{L^\infty}\|G\|^2_{L^2}\\
&\qquad+\|\mathcal{S}(\u_2)\|_{L^\infty}^2\|\mathcal{T}(\u_2)\|^2_{L^\infty}\|G\|^2_{L^2}
+\|r\|^2_{L^2}\|\mathcal{T}(\u_1)\|_{L^\infty}\|\mathcal{T}(\u_2)\|_{L^\infty}^2\\
&\le 5\varepsilon\|\nabla
v\|^2_{L^2}+\eta_3(\varepsilon)(\|r\|^2_{L^2}+\|v\|^2_{L^2}+\|G\|^2_{L^2})
\end{split}
\end{equation}
with
\begin{equation*}
\begin{split}
\eta_3(\varepsilon)&=C(\varepsilon)\|\mathcal{S}(\u_1)\|^2_{L^\infty}\|\u_1\|^2_{L^\infty}+C(\varepsilon)\|\partial_t
\u_2\|_{L^3}^2+ \|\mathcal{S}(\u_1)\|_{L^\infty}\|\nabla
\u_1\|_{L^\infty}\\&\qquad+2\|\u_2\|_{L^\infty}\|\nabla
\u_1\|_{L^\infty}
+C(\varepsilon)\|\mathcal{S}(\u_2)\|_{L^\infty}^2\|\u_2\|_{L^\infty}^2\\
&\qquad+C(\varepsilon)(\sup\{P'(\eta): C(T)^{-1}\le\eta\le
C(T)\})^2\\&\qquad+C(\varepsilon)(\|\mathcal{S}(\u_1)
\|_{L^\infty}^2\|\mathcal{T}(\u_1)\|^2_{L^\infty}+\|\mathcal{S}(\u_2)\|_{L^\infty}^2\|\mathcal{T}(\u_2)\|^2_{L^\infty})
\\&\qquad+\|\mathcal{T}(\u_1)\|_{L^\infty}^2\|\mathcal{T}(\u_2)\|_{L^\infty}^2.
\end{split}
\end{equation*}

Summing up \eqref{1213}, \eqref{1212}, and \eqref{15}, by taking
$\varepsilon=\f{\mu}{20}$, we obtain
\begin{equation}\label{16}
\begin{split}
\f{d}{dt}\int_{\R^3} &(\mathcal{S}(v_1)|v|^2+|r|^2+|G|^2)dx+\mu\int_{\R^3}|\nabla v|^2dx\\
&\le 2(\eta_3(\varepsilon)+\eta_2(\varepsilon)+\eta_1(\varepsilon))(\|v\|_{L^2}^2+\|r\|^2_{L^2}+\|G\|^2_{L^2})\\
&\le
2\eta(\varepsilon,t)\int_{\R^3}(\mathcal{S}(\u_1)|v|^2+|r|^2+|G|^2)dx,
\end{split}
\end{equation}
with
$$\eta(\varepsilon, t)=\f{\eta_3(\varepsilon)+\eta_2(\varepsilon)+\eta_1(\varepsilon)}{\min\{\min_{x\in\R^3}\mathcal{S}(\u_1)(x,t),1\}}.$$
It is a routine matter to establish the integrability with respect
to $t$ of the function $\eta(\varepsilon,t)$ on the interval
$(0,T)$. This is a consequence of the regularity of $\u_1, \u_2\in
\mathcal{W}$ and the estimates in Lemmas \ref{r} and \ref{g} for
$\mathcal{S}(\u_i)$, $\mathcal{T}(\u_i)$ with $i=1,2.$ Therefore,
\eqref{16}, combining with Gronwall's inequality, implies
$$\int_{\R^3}(\mathcal{S}(\u_1)|v|^2+|r|^2+|G|^2)dx=0,\quad \textrm{for all}\quad t\in(0,T),$$
and consequently
$$v\equiv 0,\quad r\equiv 0,\quad G\equiv 0.$$
Thus, the proof of uniqueness is complete.

\bigskip\bigskip

\section*{Acknowledgments}

Xianpeng Hu's research was supported in part by the National
Science Foundation grant DMS-0604362 and by the Mellon Predoctoral
Fellowship at the University of Pittsburgh.
Dehua Wang's research was supported in part by the National Science
Foundation under Grants DMS-0604362 and DMS-0906160, and by the Office of Naval
Research under Grant N00014-07-1-0668.

\end{document}